\documentclass{amsart}
\usepackage[T1]{fontenc}
\usepackage{textcomp}
\usepackage[utf8x]{inputenc}
\usepackage[swedish,english]{babel}
\usepackage{ucs}
\usepackage{graphicx}
\usepackage{mathrsfs}
\usepackage{rotating} 
\usepackage{tikz}
\usepackage{amssymb}
\usepackage[titletoc,toc]{appendix}
\usepackage{hyperref}
\usepackage{longtable}
\usepackage{float}
\numberwithin{equation}{section} 
\makeatletter
\def\nobreakhline{%
  \noalign{\ifnum0=`}\fi
    \penalty\@M
    \futurelet\@let@token\LT@@nobreakhline}
\def\LT@@nobreakhline{%
  \ifx\@let@token\hline
    \global\let\@gtempa\@gobble
    \gdef\LT@sep{\penalty\@M\vskip\doublerulesep}
  \else
    \global\let\@gtempa\@empty
    \gdef\LT@sep{\penalty\@M\vskip-\arrayrulewidth}
  \fi
  \ifnum0=`{\fi}%
  \multispan\LT@cols
     \unskip\leaders\hrule\@height\arrayrulewidth\hfill\cr
  \noalign{\LT@sep}%
  \multispan\LT@cols
     \unskip\leaders\hrule\@height\arrayrulewidth\hfill\cr
  \noalign{\penalty\@M}%
  \@gtempa}
\makeatother
\newtheorem{thm}{Theorem}[section]
\newtheorem{example}{Example}[section]
\newtheorem{cor}[thm]{Corollary}

\newtheorem{alg}[thm]{Algorithm}
\theoremstyle{definition}
\newtheorem{definition}[thm]{Definition}
\theoremstyle{remark}
\newtheorem{rem}[thm]{Remark}
\numberwithin{equation}{section}

\newcommand{\Hom}{\mathrm{Hom}}

\newcommand{\Tr}{\mathrm{Tr} \!}
\newcommand{\undcal}[2]{#1_{\mathcal{#2}}}
\newcommand{\intpos}[2]{\mathscr{#1}\! \left(\mathcal{#2}\right)}
\newcommand{\gintpos}[3]{\mathscr{#1}^{#2}\left(\mathcal{#3}\right)}
\newcommand{\derhamc}[2]{H^{#1}_c \! \left(#2\right)}
\newcommand{\derham}[2]{H^{#1} \! \left(#2\right)}

\newcommand{\cd}{\mathrm{cd}}

\makeindex
\begin{document}

\title[Toric arrangements associated to root systems]
{Cohomology of complements of toric arrangements associated to root systems}%
\author{Olof Bergvall}%
\address{Matematiska institutionen, Uppsala universitet, Box 480, 751 06 Uppsala, Sweden}
\email{olof.bergvall@math.uu.se}

\begin{abstract}
   We develop an algorithm for computing the cohomology of complements of toric arrangements.
   In the case a finite group $\Gamma$ is acting on the arrangement, the algorithm determines the
   cohomology groups as representations of $\Gamma$.
   As an important application, we determine the cohomology groups of the complements of
   the toric arrangements associated to root systems of exceptional type as representations
   of the corresponding Weyl groups.
   
   \vspace{5pt}
   \textit{Weyl groups}.
\end{abstract}

\maketitle


\section{Introduction}
An \emph{arrangement} is a finite set of closed subvarieties of a variety. 
Despite their simple definition, arrangements are of interest to a wide range of areas of 
mathematics such as algebraic geometry, topology, combinatorics, Lie theory and singularity theory.

Classically, most attention has been given to arrangements of hyperplanes in an affine space
but in the last two decades an increasing amount of interest has been drawn towards
\emph{toric arrangements}, i.e. arrangements $\mathcal{A}=\{A_i\}_{i \in I}$ of 
codimension one subtori inside an ambient torus $X$.
In particular, a lot of work has gone into understanding the cohomology of
complements $X_{\mathcal{A}}=X \setminus \cup_{i \in I} A_i$ of toric arrangements. 
This is also the topic of the present paper.

More specifically we are interested in the following situation. Let $\Gamma$ be a finite group
of automorphisms of $X$ which stabilize $\mathcal{A}$ as a set.
The group $\Gamma$ then acts on the complement $\undcal{X}{A}$ and this action
induces a linear action on the cohomology groups $\derham{i}{\undcal{X}{A}}$. In other words, 
the groups $\derham{i}{\undcal{X}{A}}$ become representations of $\Gamma$.
A typical example (and the example which will be given most attention in the present paper) 
is when $\Gamma$ is a Coxeter group, $L$ is a lattice with a $\Gamma$-action,
$X=L \otimes \mathbb{C}^*$ and $\mathcal{A}$ is the set of subtori fixed by reflections
of $\Gamma$. 

The main result of this work is an efficient algorithm, see Algorithms~\ref{posetalg} and \ref{cohalg},
which computes $\derham{i}{\undcal{X}{A}}$ together with its representation structure.
The algorithm is applicable rather generally but as one application we determine the cohomology
of the complements of the toric arrangements associated to root systems of exceptional type
as representations of the corresponding Weyl groups. This is not only of interest in its own
but also of interest in moduli theory. For instance:
\begin{itemize}
 \item the toric arrangement associated to $E_6$ is closely related to moduli spaces
 of cubic surfaces (see e.g. \cite{bergvallgounelas}),
 \item the toric arrangement associated to $E_7$ is closely related to moduli spaces of
 quartic curves, abelian threefolds and Del Pezzo surfaces of degree $2$ (see e.g. 
 \cite{bergstrombergvall}, \cite{bergvall}, \cite{looijenga}),
 \item the toric arrangement associated to $E_8$ is closely related to moduli spaces
 of Del Pezzo surfaces of degree $1$ and to moduli of genus $4$ curves
 (see e.g. \cite{colombovangeemenlooijenga}).
\end{itemize}
In fact, the results of the present paper has been an essential ingredient
in determining the cohomology of many of the moduli spaces mentioned above
(others, in particular the case of Del Pezzo surfaces of degree $1$, is work in progress).
We also mention that recent results suggest a connection between toric arrangements and
cluster varieties, log-Calabi-Yau varieties and mirror symmetry (see e.g. \cite{grosshackingkeel} and
\cite{looijenga}).

The study of cohomology groups of complements of
toric arrangements can be traced back all the way to Arnol'd \cite{arnold} but
the theory in its current form is probably best attributed to Lehrer \cite{lehrer}
or Looijenga \cite{looijenga}. In an abstract sense, Looijenga determined the cohomology
groups of complements of toric arrangements - from this perspective our main contribution lies
in the explicit nature of the results and the novel approach in the method.
We also mention that Looijenga's work spurred intense investigations and reformulations
from more combinatorial perspectives, starting with the work of De Concini and
Procesi \cite{deconciniprocesi2005}.

The present paper is to a large extent a toric analogue of the paper
\cite{fleischmannjaniszczak} of Fleischmann and Janiszczak
where they do the corresponding computations in the case of hyperplane arrangements.
We refer to Remark~\ref{fjrem} for a more thorough comparison with their work
but mention already now that while the hyperplane case is entirely determined
by the intersection poset, see Definition~\ref{intposdef}, we also need to take
``local'' topological and arithmetic data into account in the toric setting.
The paper \cite{eisenbudsturmfels} by Eisenbud and Sturmfels provides many
of the tools necessary to carry out these computations.

The paper is organized as follows. In Section~\ref{genarrsec} we introduce terminology
around general arrangements and review how the cohomology of the complements of a minimally pure
arrangement can be determined via inclusion-exclusion techniques.
In Section~\ref{torarrsec} we specialize to the situation of a toric arrangement.
We show that the cohomology of the complement of a toric arrangement can be computed
using only arithmetic of $\mathbb{Z}$-modules. On the basis of this we develop
Algorithms~\ref{posetalg} and \ref{cohalg} which when combined computes the
cohomology of the complement of a toric arrangement. This computation is equivariant
with respect to a finite group $\Gamma$ acting on the ambient torus and stablizing the arrangement.
Finally, in Section~\ref{appedixa} we use these algorithms to compute the cohomology groups of the complements
of the toric arrangements associated to exceptional root systems as representations
of the corresponding Weyl groups.

\section{General arrangements}
 \label{genarrsec}
 Unless otherwise specified, we shall always work over the complex numbers.

\begin{definition}
 Let $X$ be a variety. An \emph{arrangement} $\mathcal{A}$ in $X$ is  a finite set 
 $\{A_i\}_{i \in I}$ of closed subvarieties of $X$ of codimension 1.
\end{definition} 
 
Given an arrangement $\mathcal{A}$ in a variety $X$ one may define its \emph{divisor}
\begin{equation*}
 \undcal{D}{A} = \bigcup_{i \in I} A_i \subset X,
\end{equation*}
and its \emph{open complement}
\begin{equation*}
 \undcal{X}{A} = X \setminus \undcal{D}{A}.
\end{equation*}
The variety $\undcal{X}{A}$ will be our main object of study.

Many interesting properties of the variety $\undcal{X}{A}$ can be deduced
from properties of $\undcal{D}{A}$ via inclusion-exclusion arguments.
The object that governs the principle of inclusion and exclusion in this
setting is the \emph{intersection poset} of $\mathcal{A}$.

\begin{definition}
\label{intposdef}
 Let $\mathcal{A}$ be an arrangement in a variety $X$. 
 The intersection poset of $\mathcal{A}$ is the set
 \begin{equation*}
  \intpos{L}{A} = \left\{ \left. \cap_{j \in J}  A_j \right| J \subseteq I  \right\}.
 \end{equation*}
 of intersections of elements of $\mathcal{A}$, ordered by reverse inclusion.
 We include $X$ as an element of $\intpos{L}{A}$ corresponding to
 the empty intersection.
\end{definition}

\begin{rem}
 The definition of the poset $\intpos{L}{A}$ is deceivingly similar to a poset
 used in many combinatorial texts. The difference is that we allow irreducible and
 even disconnected elements in $\intpos{L}{A}$ while the ``combinatorial'' counterpart
 consists of irreducible components of intersections.
 The reason why we have not stuck to this convention will become apparent
 after Definition~\ref{modposdef}.
 
 It might also be of interest to compare our approach to the theory
 of arithmetic matroids of Br\"andén and Moci \cite{brandenmoci} (which is more 
 closely related to what is above referred to as the ``combinatorial intersection poset'')
 and to the theory of matroids over rings of Fink and Moci \cite{finkmoci}
 (which is more closely related to our approach - our approach is almost
 an embedded, dual version of a matroid over $\mathbb{Z}$). 
\end{rem}

Since $\intpos{L}{A}$ is a poset, it has a M\"obius function
$\mu: \intpos{L}{A} \times \intpos{L}{A} \to \mathbb{Z}$ defined inductively by
setting $\mu(Z,Z) = 1$ and 
\begin{equation*}
  \sum_{Y \geq Z' \geq Z} \mu(Z',Z) =0, \quad \text{if } Y \geq Z,
\end{equation*}
where the sum is over all $Z' \in \intpos{L}{A}$ between $Y$ and $Z$. 
Since we shall
exclusively be interested in the values of the M\"obius function at the minimal element $X$,
we shall use the simplified notation $\mu(Z):=\mu(Z,X)$.
See Example \ref{B2ex} for an explicit arrangement together with its
intersection poset and M\"obius function.

\subsection{Group actions on arrangements}

Let $\Gamma$ be a finite group of automorphisms of $X$ that stabilizes
$\mathcal{A}$ as a set. Such an action induces an action on $\undcal{X}{A}$
and thus yields a linear action on the de Rham cohomology groups $\derham{i}{\undcal{X}{A}}$
as well as the de Rham cohomology groups with compact supports $\derhamc{i}{\undcal{X}{A}}$.
In other words,
each cohomology group is naturally a $\Gamma$-representation.

One way to encode this information is via \emph{equivariant Poincar\'e polynomials}.
If $\Gamma$ is a finite group acting on a smooth variety $Y$ as above, we define the equivariant Poincar\'e polynomial
of $Y$ at $g \in \Gamma$ as
\begin{equation*}
 P(Y,t)(g) := \sum_{i \geq 0} \Tr \left(g, \derham{i}{Y} \right) \cdot t^i,
\end{equation*}
where $\Tr \left(g, \derham{i}{Y} \right)$ denotes the trace of
$g$ on $\derham{i}{Y}$. We define the \emph{compactly supported equivariant Poincar\'e polynomial}
$P_c(Y,t)(g)$ in a completely analogous way.

Poincar\'e polynomials are not additive so one should not expect
inclusion-\linebreak exclusion arguments to apply in the computation of 
Poincar\'e polynomials. Nonetheless, if $X$ and $\mathcal{A}$ are
nice enough this turns out to be the case. One such instance is
when both $X$ and $\mathcal{A}$ are minimally pure.
We refer to \cite{dimcalehrer} or \cite{macmeikan} for the general definitions of minimal purity
but recall that a smooth and irreducible variety $X$ is minimally pure
if $\derhamc{i}{X}$ is a pure Hodge structure of weight $2i-2\mathrm{dim}(X)$.
We also recall that the conditions for minimal purity are satisfied 
when $X$ is an affine
or projective space and each element of $\mathcal{A}$ is a hyperplane
(by results of Brieskorn \cite{brieskorn})
or when $X$ is a torus and each element of $\mathcal{A}$ is a subtorus
of codimension one (by results of Looijenga \cite{looijenga}).

For an element $g \in \Gamma$, we let $\gintpos{L}{g}{A}$ denote the subposet of $\intpos{L}{A}$
consisting of elements fixed by $g$ and we let $\mu_g$ denote the M\"obius function
of $\gintpos{L}{g}{A}$.

\begin{thm}[MacMeikan \cite{macmeikan}]
\label{macmeikanthm}
 Let $\mathcal{A}=\left\{A_i\right\}_{i \in I}$ be a minimally pure arrangement
 in a minimally pure variety $X$ and let $\Gamma$ be a finite group of automorphisms
 of $X$ that stabilizes $\mathcal{A}$ as a set. Then, for each $g \in \Gamma$
 \begin{equation*}
  P_c(\undcal{X}{A},t)(g) = \sum_{Z \in \gintpos{L}{g}{A}} \mu_g(Z) (-t)^{\cd(Z)} P_c(Z,t)(g),
 \end{equation*}
 where $\cd(Z)$ denotes the codimension of $Z$ in $X$.
\end{thm}

We will use the following version of the above result.

\begin{cor}
\label{macmeikancor}
 Let $\mathcal{A}=\left\{A_i\right\}_{i \in I}$ be a minimally pure arrangement
 in a minimally pure variety $X$ and let $\Gamma$ be a finite group of automorphisms
 of $X$ that stabilizes $\mathcal{A}$ as a set. Suppose also that both $\undcal{X}{A}$
 and each element of $\intpos{L}{A}$ satisfy Poincar\'e duality. Then, for
 each $g\in \Gamma$
 \begin{equation*}
  P(\undcal{X}{A},t)(g) = \sum_{Z \in \gintpos{L}{g}{A}} \mu_g(Z) (-t)^{\cd(Z)} P(Z,t)(g),
 \end{equation*}
 where $\cd(Z)$ denotes the codimension of $Z$ in $X$.
\end{cor}

\begin{proof}
 Poincar\'e duality tells us that if $M$ is a smooth manifold of complex dimension
 $n$, then $H^k(M) = H_c^{2n-k}(M)$ (see \cite{madsentornehave}). Phrased in terms of
 Poincar\'e polynomials, this translates to
 \begin{equation*}
  P_c(M,t)=t^{2n} \cdot P(M,t^{-1}).
 \end{equation*}
 We apply Poincar\'e duality to Theorem~\ref{macmeikanthm} and get
 \begin{equation*}
  t^{2n} P(\undcal{X}{A},t^{-1})(g) = \sum_{Z \in \gintpos{L}{g}{A}} \mu_g(Z) (-t)^{\cd(Z)} \cdot t^{2 \mathrm{dim}(Z)} \cdot P(Z,t^{-1})(g). 
 \end{equation*}
 We thus have that
 \begin{align*}
  P(\undcal{X}{A},t^{-1})(g) & = \sum_{Z \in \gintpos{L}{g}{A}} \mu_g(Z) (-t)^{\cd(Z)} \cdot t^{2 \mathrm{dim}(Z)-2n} \cdot P(Z,t^{-1})(g) = \\
   & = \sum_{Z \in \gintpos{L}{g}{A}} \mu_g(Z) (-t)^{\cd(Z)} \cdot t^{-2 \cd(Z)} \cdot P(Z,t^{-1})(g) = \\
   & = \sum_{Z \in \gintpos{L}{g}{A}} \mu_g(Z) (-t^{-1})^{\cd(Z)} \cdot P(Z,t^{-1})(g).
 \end{align*}
 We now arrive at the desired formula by substituting $t^{-1}$ for $t$.
\end{proof}

\begin{example}
\label{B2ex}
Let $X= \left( \mathbb{C}^* \right)^2$ and consider the arrangement $\mathcal{A}$ in $X$
consisting of the four subtori given by the equations
\begin{equation*}
 A_1: \, z_1=1, \quad A_2: \, z_1^2z_2=1, \quad A_3: \, z_1z_2^2=1, \quad A_4: \, z_2=1.
\end{equation*}
Let $\xi$ be a primitive third root of unity. We then have
\begin{align*}
 & A_1 \cap A_2 = A_1 \cap A_4 = A_3 \cap A_4 = \{(1,1)\}, \\
 & A_1 \cap A_3 = \{(1,1), \, (1,-1)\}, \\
 & A_2 \cap A_3 = \{(1,1), \, (\xi,\xi), \, (\xi^2,\xi^2)\}, \\
 & A_2 \cap A_4 = \{(1,1), \, (-1,1)\},
\end{align*}
and all further intersections are equal to $\{(1,1)\}$. We thus have
the poset $\intpos{L}{A}$ (the numbers in the upper left corners are the
values of the M\"obius function):
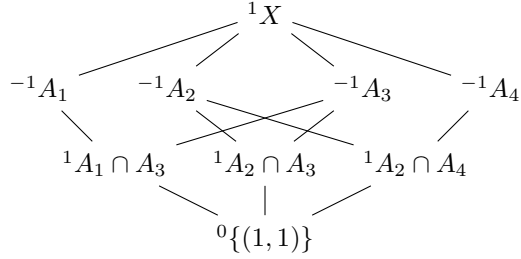
\begin{center}
\begin{tikzpicture}
    \node (X) at (0,-3) {$^1X$};
    \node (A1) at (-3,-2) {$^{-1}A_{1}$};
    \node (A2) at (-1.3,-2) {$^{-1}A_{2}$};
    \node (A3) at (1.3,-2) {$^{-1}A_{3}$};
    \node (A4) at (3,-2) {$^{-1}A_{4}$};
    \node (A13) at (-2,-1) {$^1A_{1} \cap A_3$};
    \node (A23) at (0,-1) {$^1A_{2} \cap A_3$};
    \node (A24) at (2,-1) {$^1A_{2} \cap A_4$};
    \node (A1234) at (0,0) {$^0\{(1,1)\}$};
    \draw (X) -- (A1) -- (A13) -- (A1234);
    \draw (X) -- (A2) -- (A23) -- (A1234);
    \draw (X) -- (A4) -- (A24) -- (A1234);
    \draw (X) -- (A3);
    \draw (A3) -- (A13);
    \draw (A3) -- (A23);
    \draw (A2) -- (A24);
\end{tikzpicture} 
\end{center}
The Poincar\'e polynomial of $X$ is $(1+t)^2$ and the Poincar\'e polynomial
of $A_i$ is $1+t$.
By Corollary~\ref{macmeikancor} we now get that the Poincaré polynomial of $\undcal{X}{A}$ is
\begin{align*}
 P(\undcal{X}{A},t) & = (1+t)^2+4 \cdot (-1) \cdot(-t)\cdot (1+t) +  (-t)^2 \cdot (2 +3+2)  = \\
 & = 1 + 6\, t + 12 \, t^2. 
\end{align*}
\end{example}

\section{Toric arrangements}
\label{torarrsec}

\begin{definition}
 Let $X$ be an $n$-torus. An arrangement $\mathcal{A}$ in $X$
 is called a \emph{toric arrangement} if each element of $\mathcal{A}$
 is a hypertorus, i.e. a subtorus of codimension one.
\end{definition}

More explicitly, a toric arrangement in $X=(\mathbb{C}^*)^n$ is determined
by a finite set $\{\chi_i\}_{i \in I}$ of characters $\chi_i$ in the character
lattice $\Hom(X,\mathbb{C}^*)$ of $X$. The subtorus $A_i$ is then given as
$A_i=\chi_i^{-1}(1)$.

For our purposes it is more convenient to start the other way around, i.e.
to start with a free $\mathbb{Z}$-module $M$ of rank $n$ and construct
the torus $X$ as $X= \Hom(M,\mathbb{C}^*)$. An arrangement in $X$ is then
given by a finite set $R=\{m_i\}_{i \in I}$ of elements $m_i \in M$ and
the subtorus $A_i$ is given as the kernel of the evaluation map at $m_i$, i.e.
\begin{equation*}
 A_i = \{\chi \in X | \chi(m_i)=1\}. 
\end{equation*}
Let $B=\{b_1, \ldots, b_n\}$ be an ordered set of free generators of $M$. An element
$\chi \in X$ is then entirely determined by its values at $b_1, \ldots, b_n$
and the coordinate ring of $X$ is therefore the ring of Laurent polynomials
\begin{equation*}
 \mathbb{C}[X] = \mathbb{C}[z_1^{\pm 1}, \ldots, z_n^{\pm 1}]
\end{equation*}
where the variable $z_i$ corresponds to $b_i$. If $m$ is an element in $R$, then
there are integers $r_1, \ldots, r_n$ such that
\begin{equation*}
 m = r_1 b_1 + \cdots + r_n b_n
\end{equation*}
and we have $\chi(m)=1$ if and only if
\begin{equation*}
 z_1^{r_1} \cdots z_n^{r_n} = 1.
\end{equation*}
We denote the Laurent polynomial $ z_1^{r_1} \cdots z_n^{r_n} - 1$ by $f_m$
and the corresponding subtorus in $\mathcal{A}$ is denoted both as $A_m$ and
$V(f_m)$. We will also make the identification $A_{m_i}=A_i$ whenever convenient.

Let $Z$ be an element of $\intpos{L}{A}$. Then $Z$ is an intersection
\begin{equation*}
 Z = \bigcap_{m \in S} A_m = \bigcap_{m \in S} V(f_m)
\end{equation*}
where $S$ is a subset of $R$. We define the ideal
\begin{equation*}
 I_{S} = \left( f_m \right)_{m \in S} \subseteq \mathbb{C}[X]
\end{equation*}
so that $Z=V(I_{S})$. The ideal $I_S$ is entirely
determined by the exponents occurring in the various Laurent polynomials $f_{\alpha}$ 
generating it or, in other words, the coefficients occurring in the elements
of $S$ when expressed in terms of the generators in $B$. Thus, if
we define the \emph{module of exponents}
\begin{equation*}
 N_S := \mathbb{Z}\left\langle S \right\rangle \subseteq M,
\end{equation*}
then the module $N_S$ determines $I_S$ and
\begin{equation*}
 V(I_S) = \mathrm{Hom}(M/N_S, \mathbb{C}^*) \subseteq \mathrm{Hom}(M, \mathbb{C}^*) =T.
\end{equation*}
For more details, see \cite{eisenbudsturmfels}.

From the above we see that there is an intimate relationship between
elements of $\intpos{L}{A}$ and submodules of $M$ generated by
subsets of $R$. We have that an inclusion $N_S \subset N_{S'}$ gives a surjection
$M/N_S \twoheadrightarrow M/N_{S'}$ and thus an inclusion
$\mathrm{Hom}(M/N_{S'},\mathbb{C}^*) \hookrightarrow \mathrm{Hom}(M/N_S,\mathbb{C}^*)$.

\begin{definition}
\label{modposdef}
 Let $R$ be a finite set of nonzero elements of a free $\mathbb{Z}$-module $M$. 
 The poset of modules of exponents is the set
 \begin{equation*}
  \mathscr{P}(R) = \left\{\mathbb{Z}\langle S \rangle | S \subseteq R \right\},
 \end{equation*}
 ordered by inclusion. 
\end{definition}

By construction we have that the poset $\mathscr{P}(R)$ is isomorphic to
the poset $\intpos{L}{A}$.
One benefit of considering the posets $\mathscr{P}(R)$ instead is that they are
more easily computed. Another benefit is that the elements
of $\mathscr{P}(R)$ are directly related to the cohomology groups of the 
corresponding elements of $\intpos{L}{A}$. In order to explain how, we recall
some basic facts about cohomology of tori that can be found in most texts on
toric varieties, e.g. Chapter 9 of \cite{coxlittleschenck}.

Let $L$ be a free $\mathbb{Z}$-module. The torus 
$T_L=\Hom(L,\mathbb{C}^*)$ has cohomology
given by
\begin{equation*}
 \derham{i}{T_L}= \bigwedge^i \derham{1}{T_L} = \bigwedge^i L.
\end{equation*}
Suppose $L'$ is another free $\mathbb{Z}$-module and let 
$T_{L'}= \Hom(L',\mathbb{C}^*)$. There is a one-to-one correspondence between homomorphisms
$L \to L'$ of free $\mathbb{Z}$-modules and morphisms 
$T_{L'} \to T_L$ and the induced
map $\derham{i}{T_L} \to \derham{i}{T_{L'}}$ is the map
\begin{equation*}
 \bigwedge^i L \to \bigwedge^i L'.
\end{equation*}
A module $M/N_S$ will not always be free but will still determine
the cohomology of $V(I_S)$ in a sense very similar to the above.
Let $L$ be a free $\mathbb{Z}$-module, $N \subset L$ a submodule and let $Q=L/N$.
The module $Q$ will split as a direct sum $Q=Q^T \oplus Q^F$,
where $Q^T$ is the torsion part and $Q^F$ is the free part of $Q$.
The variety $T_Q=\Hom(Q, \mathbb{C}^*)$ consists of $\left| Q^T \right|$ connected components, each
isomorphic to $\Hom \! \left(Q^F, \mathbb{C}^* \right)$. 
The $i$th cohomology group of $T_Q$ is given by
\begin{equation*}
 \derham{i}{T_Q} = \bigoplus_{v \in Q^T} \bigwedge^i Q^F.
\end{equation*}
Let $\varphi: L \to L'$ be a homomorphism and define $Q'=L'/\varphi(N)$.
The morphism $\varphi$ induces a morphism $Q \to Q'$ which in turn gives rise to a morphism $T_{Q'} \to T_Q$
and the induced map $\derham{i}{T_Q} \to \derham{i}{T_{Q'}}$ is the map
\begin{equation*}
\bigoplus_{v \in Q^T} \bigwedge^i Q^F \to  \bigoplus_{v' \in Q'^T} \bigwedge^i Q'^F.
\end{equation*}

In conclusion, in order to use Corollary~\ref{macmeikancor} to compute the cohomology
of $\undcal{X}{A}$ we can first compute the poset $\mathscr{P}(R)$ of submodules of $M$ generated
by subsets of $R$ and then for each element $N_S$ of $\mathscr{P}(R)$ compute 
the torsion part of $Q_S=M/N_S$ and all wedge powers of the free part of $Q_S$. 
We will elaborate on this in the following while also introducing group actions.

\subsection{Group actions on toric arrangements}

Let $\Gamma$ be a group acting linearly on $M$ (from the left) such that
$\Gamma$ fixes $R$ as a set. Then $\Gamma$ will also act on $\undcal{X}{A}$
(from the right) by precomposition, i.e.
\begin{equation*}
(\chi . g)(v) =  \chi(g.v).
\end{equation*}
We write $\mathscr{P}^g(R)$ to denote the subposet of $\mathscr{P}(R)$ of
modules fixed by $g$. Just as $\mathscr{P}(R)$ is isomorphic to $\intpos{L}{A}$
we have that $\mathscr{P}^g(R)$ is isomorphic to $\gintpos{L}{g}{A}$.

Let $g$ be an element of $\Gamma$. If $N$ is an element
of $\mathscr{P}^g(R)$, then $N \cap R$ is a union of $g$-orbits of $R$.
Since $N=\mathbb{Z}\langle N \cap R \rangle$,
we may compute $\mathscr{P}^g(R)$ via the following steps.

\begin{alg}
\label{posetalg}
Let $R$ and $\Gamma$ be as above and let $g$ be an element
of $\Gamma$. Then the following algorithm computes the poset $\mathscr{P}^g(R)$:
\begin{itemize}
 \item[(1)] Compute the $g$-orbits of $R$.
 \item[(2)] Compute the set $\mathscr{P}^g(R)_{\mathrm{set}}$ of all (distinct) $\mathbb{Z}$-spans of unions of $g$-orbits.
 \item[(3)] Investigate the inclusion relations of the elements of $\mathscr{P}^g(R)_{\mathrm{set}}$.
\end{itemize}
\end{alg}

In practice, we represent elements of $\mathscr{P}^g(R)_{\mathrm{set}}$ by their
Hermite normal forms (i.e. the nontrivial part of the echelon form of any matrix with a spanning 
set as rows which need to be computed with some care since we are working over the integers).
Thus, step (2) consists of performing Gaussian elimination on matrices and making
sure we only save each matrix once. Step (3) is by far the most computationally demanding since,
in principle, one needs to make $|\mathscr{P}^g(R)|^2$ checks of whether one module
is a submodule of another.
Fortunately, it is very parallelisable.

\begin{rem}
 Algorithm~\ref{posetalg} is a toric analogue of the algorithm of Fleischmann and Janiszczak
 \cite{fleischmannjaniszczak} for computing the intersection poset associated to a hyperplane
 arrangement.
\end{rem}

The group  $\Gamma$ acts not only on poset $\mathscr{P}(R)$ but also 
on the cohomology groups of the individual elements of $\mathscr{P}(R)$.
Let $g$ be an element of $\Gamma$ and let $N$ be an element of $\mathscr{P}^g(R)$.
Let $Q_N=M/N$ and let $Q^T_N$ and $Q^F_N$ denote the torsion part and free part, respectively.
Recall that the $i$th cohomology group of the element $Z_N\in \gintpos{L}{g}{A}$ corresponding to $N$
is given by
\begin{equation*}
 H^i(Z_N) = \bigoplus_{v \in Q^T_N} \bigwedge^i Q_N^F.
\end{equation*}
The element $g$ may both permute the elements of $Q^T_N$ (i.e. permute
the components of $Z_N$) and act on $\wedge^i Q_N^F$ via its action on
$Q^F_N$.

More explicitly, the module $Q_N$ is determined by the echelon basis matrix of $N$.
Torsion elements of $Q_N$ stems from rows in the echelon basis matrix whose entries has
a greatest common divisor greater than $1$. The module $Q_N^F$ is the module 
$M/\mathrm{Sat}(N)$,
where $\mathrm{Sat}(N)=N \otimes_{\mathbb{Z}} \mathbb{Q} \cap M$ is the saturation of $N$, i.e. the module generated by the elements obtained from the Hermite normal form of $N$
by dividing each row by the greatest common divisor of its entries (the division procedure does not
immediately yield a Hermite normal form of $\mathrm{Sat}(N)$).
A row $(r_1, \ldots, r_n)$ in the echelon basis matrix of $N$ corresponds to the equation
\begin{equation*}
 z_1^{r_1} \cdots z_n^{r_n} = 1.
\end{equation*}
If $\mathrm{gcd}(r_1, \ldots, r_n)=d$ we may write $r_i = d \cdot r'_i$ and
\begin{equation*}
 (z_1^{r'_1} \cdots z_n^{r'_n})^d = 1,
\end{equation*}
i.e. an equation for $d$ non-intersecting hypertori, namely the hypertorus
given by
\begin{equation*}
z_1^{m'_1} \cdots z_n^{m'_n}=1, 
\end{equation*}
translated by multiplication by powers
of a primitive $d$th root of unity.

A linear map $g: M \to M$ which fixes $N$ can be analysed in two steps.
Firstly, we can investigate how it ``permutes different roots of unity'',
more precisely, how it acts on $Q_N^T$. The elements of $Q_N^T$ correspond
to connected components of $Z_N$ and a component is fixed by $g$ if and
only if the corresponding element of $Q_N^T$ is fixed. Of course,
only fixed components can contribute to the trace of $g$ on $\derham{i}{Z_N}$.
Once we have determined which of the components that are fixed it
suffices to compute the trace of $g$ on the cohomology on one of those components,
e.g. the component corresponding to the zero element of $Q_N^T$.

We may now write down an algorithm for computing the equivariant Poincar\'e polynomial of
an element $Z \in \intpos{L}{A}$.

\begin{alg}
\label{cohalg}
Let $\Gamma$, $M$ and $R$  be as above and let $Z \in \intpos{L}{A}$ correspond
 to the submodule $N$ in $\mathscr{P}(R)$. Let $g$ be an element of $\Gamma$ stabilizing $Z$.
 Then $P(Z,t)(g)$ can computed via the following steps.
 \begin{itemize}
  \item[(1)] Compute the number $m$ of elements in $Q_N^T$ which are fixed by $g$
  (for instance by lifting each element of $Q_N^T$ to $M$, acting on the lifted element by $g$ and pushing the result down to $Q_N^T$).
  \item[(2)] Compute $\mathrm{Tr}(g,Q_N^F)$ (for instance as $\mathrm{Tr}(g,M)-\mathrm{Tr}(g,N)$).
  \item[(3)] Using the knowledge of $\mathrm{Tr}(g,Q_N^F)$, compute $\mathrm{Tr}(g,\wedge^i Q_N^F)$ for $i=1, \ldots, n$ (for instance via
  the Newton-Girard method, i.e. iteratively using Newton's identities).
 \end{itemize}
 The polynomial $P(Z,t)(g)$ is now given by
 \begin{equation*}
  P(Z,t)(g) = m \cdot \sum_{i=0}^n \mathrm{Tr}(g,\wedge^i Q_N^F) t^i.
 \end{equation*}
\end{alg}

 \section{Toric arrangements associated to exceptional root systems}
 \label{appedixa}
 We now explain how to associate a toric arrangement to a root system.
 Let $\Phi$ be a root system and let $M=\mathbb{Z}\langle \Phi \rangle$ be the 
 $\mathbb{Z}$-linear span of the elements of $\Phi$, i.e. the root lattice. Let $X$ be the torus
 $X=\mathrm{Hom}(M,\mathbb{C}^*)$ and, for each $\alpha \in \Phi$ let
 $A_{\alpha}$ be the kernel of the evaluation map $X \to \mathbb{C}^*$ defined
 by $\alpha$, i.e.
 \begin{equation*}
  A_{\alpha} = \{\chi \in X| \chi(\alpha) = 1\}.
 \end{equation*}
 We note that if $\alpha'=-\alpha$, then $A_{\alpha'}=A_{\alpha}$.
 The $A_{\alpha}$ are codimension $1$ subtori and constitute a toric arrangement 
 $\mathcal{A}_{\Phi} = \{A_{\alpha}\}_{\alpha \in \Phi^{+}}$, where
 $\Phi^{+}$ is a choice of positive roots for $\Phi$. For notational convenience,
 we denote the complement $X_{\mathcal{A}_{\Phi}}$ by $X_{\Phi}$.
 
 We have used Algorithms~\ref{posetalg} and \ref{cohalg}
 to construct a \texttt{SageMath} \cite{sagemath} program computing the cohomology
 of the complement of a toric arrangement associated to a root system as a representation
 of the corresponding Weyl group. The code can be found at the permanent repository
 \url{https://github.com/OlofBergvall/CohTorArr}. 
 The validity of this program has been tested against the results in
 \cite{bergvall_an}, \cite{ardilacastillohenley}, \cite{lehrer_toral} and \cite{liutranyoshinaga}
 as well as several other works which are referred to more properly below.
 In this section we present the results produced by this program for the exceptional root systems.
 
 The irreducible representations of Weyl groups of exceptional root systems
 can be described in terms of two integers, $d$ and $e$. The integer $d$ is
 the degree of the representation. The integer $e$ can be defined using
 the standard representation $\chi_{\mathrm{std}}$ as follows, see \cite{carter}, p. 411.
 
 \begin{definition}
  Let $\Phi$ be an exceptional root system and let $W$ be its Weyl group.
  Then each irreducible representation of $W$ occurs in some symmetric power $\mathrm{Sym}^i\chi_{\mathrm{std}}$
  of $\chi_{\mathrm{std}}$. Given an irreducible representation $\chi$, let $e$
  be the integer such that $\chi$ occurs as a direct summand in $\mathrm{Sym}^e\chi_{\mathrm{std}}$
  but not in $\mathrm{Sym}^i\chi_{\mathrm{std}}$ for all $i<e$.
 \end{definition}

 For root systems of type $E$, the integers $d$ and $e$ uniquely determine the
 irreducible representations. However, for the root systems $F_4$ and $G_2$ this
 is not the case.
 We denote the irreducible character corresponding
 to $d$ and $e$ by $\phi_d^e$.
 If there are two characters corresponding
 to the same $d$ and $e$ we add a second subscript to distinguish between them.
 
 We refer to \cite{carter} for character tables of the Weyl groups
 of type $G_2$ and $F_4$. It should be noted that Carter denotes the
 characters $\phi_d^e$, $\phi_{d,1}^e$ and $\phi_{d,2}^e$ by
 $\phi_{d,e}$, $\phi_{d,e}'$ and $\phi_{d,e}''$. We have chosen
 to denote the characters differently because we are in need of notational compactness.

 \begin{rem}
 \label{fjrem}
  This work parallels that of Fleischmann and Janiszczak \cite{fleischmannjaniszczak} in the hyperplane case
  and we therefore comment on a few differences. 
  Most obviously, in the hyperplane case the elements of the intersection poset
  are linear subspaces and therefore have trivial cohomology. The cohomology is therefore
  determined entirely by the combinatorics of the poset. This is not the case in the toric setting.
  Furthermore, in the hyperplane case one can reduce modulo ``good'' primes and do computations
  over finite fields. This speeds up computations dramatically. In the toric setting we need
  to do all computations over $\mathbb{Z}$.
  Finally, one may observe that the poset $\mathscr{P}^g(R)$ is largest when $g$ is the identity
  element or minus the identity (more complicated group elements typically have fewer fixed elements).
  In the hyperplane case, the Poincaré polynomials $P(X_{\Phi},t)(\pm \mathrm{id})$ were known
  (and known to be the same). Thus, the computationally most involved cases were already taken care of
  by other methods. In the toric case we only knew $P(X_{\Phi},t)(\mathrm{id})$ from before which
  forces us to consider the most complicated posets. This increases the computational complexity more than 
  10 000 times. On the other hand, there has also been a significant increase in computational power
  over the two decades separating our work from that of Fleischmann and Janiszczak.
 \end{rem}
 

 \subsection{The root system $G_2$}
 Let $\Phi$ be the root system of type $G_2$. Then the Poincar\'e polynomial
 of $X_{\Phi}$ is
 \begin{equation*}
  P(X_{\Phi},t) = 19 \, t^2 + 8 \, t + 1.
 \end{equation*}
 The cohomology of $X_{\Phi}$ as a representation of the Weyl group
 of $G_2$ is given in Table~\ref{g2cohtable}.
 
 \begin{center}
 \begin{table}[H]
  \label{g2cohtable}
 \caption{The cohomology groups of the complement of the toric arrangement associated to $G_2$ as
 representations of the Weyl group.}
 \begin{tabular}{rrrrrrr}
  \, & $\phi_{1}^0$ & $\phi_{1}^6$ & $\phi_{1,1}^3$ & $\phi_{1,2}^3$ & $\phi_{2}^1$ & $\phi_{2}^2$ \\
  \hline
 $H^0$ & 1 &$\cdot$&$\cdot$&$\cdot$&$\cdot$&$\cdot$\\
 $H^1$ & 2 &$\cdot$&$\cdot$&$\cdot$& 1 & 2 \\
 $H^2$ & 2 & 1 & 1 & 1 & 3 & 4
 \end{tabular}
 \end{table}
\end{center}

 \subsection{The root system $F_4$}
 Let $\Phi$ be the root system of type $F_4$. Then the Poincar\'e polynomial
 of $X_{\Phi}$ is
 \begin{equation*}
  P(X_{\Phi},t) = 2153 \, t^{4} + 1260 \, t^{3} + 286 \, t^{2} + 28 \, t + 1 .
 \end{equation*}
 The cohomology of $X_{\Phi}$ as a representation of the Weyl group
 of $G_2$ is given in Table~\ref{f4cohtable}.
 
  \begin{rem}
  The polynomial $P(X_{\Phi},t)$ was computed by Moci in \cite{moci}
  as the Poincar\'e polynomial of the complement of the toric arrangement associated to the
  \emph{coroot} system of $F_4$.
 \end{rem} 
 
 \begin{center}
 \begin{table}[H]
  \label{f4cohtable}
 \caption{The cohomology groups of the complement of the toric arrangement associated to $F_4$ as
 representations of the Weyl group.}
 \begin{tabular}{rrrrrrrrrrrrrrrrrrrrrrrrrr}
 \, & $\phi_{1}^0$ & $\phi_{1}^{24}$ & $\phi_{1,1}^{12}$ & $\phi_{1,2}^{12}$ & $\phi_{2,1}^{16}$ & $\phi_{2,2}^{4}$ & $\phi_{2,2}^{16}$ & $\phi_{2,4}^{4}$ &
 $\phi_{4}^{1}$ & $\phi_{4,2}^{7}$ & $\phi_{4,1}^7$ & $\phi_{4}^{13}$ & $\phi_{4}^8$ \\
 \hline
 $H^0$ & 1 &$\cdot$&$\cdot$&$\cdot$&$\cdot$&$\cdot$&$\cdot$&$\cdot$&$\cdot$&$\cdot$&$\cdot$&$\cdot$&$\cdot$\\
 $H^1$ & 2 &$\cdot$&$\cdot$&$\cdot$&$\cdot$& 1 &$\cdot$& 1 & 1 &$\cdot$&$\cdot$&$\cdot$&$\cdot$\\
 $H^2$ & 2 &$\cdot$&$\cdot$&$\cdot$&$\cdot$& 3 &$\cdot$& 3 & 2 &$\cdot$&$\cdot$&$\cdot$& 2 \\
 $H^3$ & 3 &$\cdot$& 1 & 2 & 1 & 6 & 2 & 5 & 3 & 3 & 2 & 1 & 8 \\
 $H^4$ & 3 & 1 & 2 & 3 & 3 & 6 & 4 & 5 & 6 & 7 & 6 & 5 & 10 \\
 \hline
 \, & $\phi_{6,2}^6$ & $\phi_{6,1}^6$ & $\phi_{8,1}^{3}$ & $\phi_{8,2}^9$ & $\phi_{8,2}^3$ & $\phi_{8,1}^9$ & $\phi_{9}^{10}$ & $\phi_{9,1}^6$ & $\phi_{9,2}^6$ & $\phi_{9,}^2$ & $\phi_{12}^4$ & $\phi_{16}^5$ \\
  \hline
 $H^0$ &$\cdot$&$\cdot$&$\cdot$&$\cdot$&$\cdot$&$\cdot$&$\cdot$&$\cdot$&$\cdot$&$\cdot$&$\cdot$&$\cdot$\\
 $H^1$ &$\cdot$&$\cdot$&$\cdot$&$\cdot$&$\cdot$&$\cdot$&$\cdot$&$\cdot$&$\cdot$& 2 &$\cdot$&$\cdot$\\
 $H^2$ & 1 & 1 & 2 &$\cdot$& 2 &$\cdot$&$\cdot$& 3 & 4 & 9 & 3 & 2 \\
 $H^3$ & 7 & 9 & 7 & 5 & 8 & 4 & 7 & 12 & 15 & 20 & 16 & 12 \\
 $H^4$ & 12 & 14 & 13 & 13 & 14 & 12 & 16 & 18 & 20 & 22 & 25 & 26 
 \end{tabular}
 \end{table}
\end{center}  
 
 \subsection{The root system $E_6$}
 Let $\Phi$ be the root system of type $E_6$. Then the Poincar\'e polynomial
 of $X_{\Phi}$ is
 \begin{equation*}
  P(X_{\Phi},t) = 58555 \, t^6 + 63378 \, t^5 + 27459 \, t^4 + 6020 \, t^3 + 705 \, t^2 + 42 \, t + 1.
 \end{equation*}
 The cohomology of $X_{\Phi}$ as a representation of the Weyl group
 of $E_6$ is given in Table~\ref{e6cohtable}.
 
 \begin{rem}
  The column $\phi_{1}^0$ gives the cohomology of $X_{\Phi}/W$. It
  was first computed by Looijenga in \cite{looijenga}.
 \end{rem}
 
 \begin{center}
 \begin{table}[H]
  \label{e6cohtable}
 \caption{The cohomology groups of the complement of the toric arrangement associated to $E_6$ as
 representations of the Weyl group.}
 \begin{tabular}{rrrrrrrrrrrrrr}
 \, & $\phi_{1}^0$ & $\phi_{1}^{36}$ & $\phi_{6}^{25}$ & $\phi_{6}^1$ & $\phi_{10}^{9}$ & $\phi_{15}^{17}$ & $\phi_{15}^{16}$ & $\phi_{15}^{5}$ & $\phi_{15}^{4}$ & $\phi_{20}^{20}$ & $\phi_{20}^2$ & $\phi_{20}^{10}$ & $\phi_{24}^{12}$  \\
  \hline
$H^0$ & 1 &$\cdot$&$\cdot$&$\cdot$&$\cdot$&$\cdot$&$\cdot$&$\cdot$&$\cdot$&$\cdot$&$\cdot$&$\cdot$&$\cdot$\\
$H^1$ & 1 &$\cdot$&$\cdot$& 1 &$\cdot$&$\cdot$&$\cdot$&$\cdot$& 1 &$\cdot$& 1 &$\cdot$&$\cdot$\\
$H^2$ &$\cdot$&$\cdot$&$\cdot$& 1 &$\cdot$&$\cdot$&$\cdot$& 1 & 2 &$\cdot$& 3 &$\cdot$&$\cdot$\\
$H^3$ &$\cdot$&$\cdot$&$\cdot$& 1 & 1 &$\cdot$&$\cdot$& 5 & 4 &$\cdot$& 6 & 2 & 1 \\
$H^4$ &$\cdot$&$\cdot$&$\cdot$& 5 & 7 & 4 & 3 & 12 & 10 & 2 & 15 & 11 & 10 \\
$H^5$ & 1 &$\cdot$& 4 & 10 & 12 & 15 & 14 & 20 & 19 & 17 & 30 & 24 & 26 \\
$H^6$ & 2 & 1 & 6 & 8 & 12 & 16 & 18 & 17 & 19 & 21 & 25 & 23 & 27 \\
\hline
\, & $\phi_{24}^{6}$ & $\phi_{30}^{15}$ & $\phi_{30}^{3}$ & $\phi_{60}^{11}$ & $\phi_{60}^{5}$ & $\phi_{60}^{8}$ & $\phi_{64}^{13}$ & $\phi_{64}^{4}$ & $\phi_{80}^{7}$ & $\phi_{81}^{6}$ & $\phi_{81}^{10}$ & $\phi_{90}^{8}$ \\
\hline
$H^0$ &$\cdot$&$\cdot$&$\cdot$&$\cdot$&$\cdot$&$\cdot$&$\cdot$&$\cdot$&$\cdot$&$\cdot$&$\cdot$&$\cdot$\\
$H^1$ &$\cdot$&$\cdot$&$\cdot$&$\cdot$&$\cdot$&$\cdot$&$\cdot$&$\cdot$&$\cdot$&$\cdot$&$\cdot$&$\cdot$\\
$H^2$ &$\cdot$&$\cdot$& 2 &$\cdot$& 2 & 1 &$\cdot$& 3 &$\cdot$& 2 &$\cdot$&$\cdot$\\
$H^3$ & 4 &$\cdot$& 9 & 3 & 11 & 6 & 1 & 14 & 9 & 14 & 5 & 10 \\
$H^4$ & 16 & 9 & 24 & 23 & 38 & 30 & 21 & 45 & 45 & 50 & 36 & 50 \\
$H^5$ & 31 & 30 & 41 & 68 & 80 & 74 & 69 & 88 & 99 & 103 & 94 & 111 \\
$H^6$ & 29 & 33 & 36 & 66 & 69 & 69 & 69 & 74 & 91 & 92 & 90 & 101 \\
 \end{tabular}
 \end{table}
\end{center}

  \subsection{The root system $E_7$}
 Let $\Phi$ be the root system of type $E_7$. Then the Poincar\'e polynomial
 of $X_{\Phi}$ is
 \begin{equation*}
 \begin{array}{rl}
   P(X_{\Phi},t)= &  3842020 \, t^7 + 3479734 \, t^6 + 1328670 \, t^5 + 268289 \, t^4 +\\
  \, & + 30800 \, t^3 + 2016 \, t^2 + 70 \, t + 1.
 \end{array}
 \end{equation*}
 The cohomology of $X_{\Phi}$ as a representation of the Weyl group
 of $E_7$ is given in Table~\ref{e7cohtable}.
 
 \begin{rem}
  The column $\phi_{1}^0$ gives the cohomology of $X_{\Phi}/W$. It
  was first computed by Looijenga in \cite{looijenga} and later corrected
  by Getzler and Looijenga in \cite{getzlerlooijenga}.
 \end{rem}
 \clearpage

\begin{center}
 \begin{table}[H]
  \label{e7cohtable}
 \caption{The cohomology groups of the complement of the toric arrangement associated to $E_7$ as
 representations of the Weyl group.}
 \begin{tabular}{rrrrrrrrrrrrr}
  \, & $ \phi_{1}^0$ & $\phi_{1}^{63}$ & $\phi_{7}^{46}$ & $\phi_{7}^{1}$ & $\phi_{15}^{7}$ & $\phi_{15}^{28}$ & $\phi_{21}^{36}$ & $\phi_{21}^{3}$ & $\phi_{21}^{33}$ & $\phi_{21}^{6}$ & $\phi_{27}^{2}$ & $\phi_{27}^{37}$ \\
  \hline
  $H^0$ & 1 &$\cdot$&$\cdot$&$\cdot$&$\cdot$&$\cdot$&$\cdot$&$\cdot$&$\cdot$&$\cdot$&$\cdot$&$\cdot$\\
  $H^1$ & 1 &$\cdot$&$\cdot$& 1 &$\cdot$&$\cdot$&$\cdot$&$\cdot$&$\cdot$&$\cdot$& 1 &$\cdot$\\
  $H^2$ &$\cdot$&$\cdot$&$\cdot$& 1 &$\cdot$&$\cdot$&$\cdot$& 1 &$\cdot$& 1 & 2 &$\cdot$\\
  $H^3$ &$\cdot$&$\cdot$&$\cdot$&$\cdot$& 1 &$\cdot$&$\cdot$& 2 &$\cdot$& 3 & 3 &$\cdot$\\
  $H^4$ &$\cdot$&$\cdot$&$\cdot$& 1 & 4 &$\cdot$&$\cdot$& 3 &$\cdot$& 7 & 8 &$\cdot$\\
  $H^5$ &$\cdot$&$\cdot$&$\cdot$& 3 & 7 & 3 & 2 & 11 & 2 & 17 & 25 & 1 \\
  $H^6$ & 2 &$\cdot$& 4 & 9 & 14 & 18 & 19 & 26 & 16 & 34 & 50 & 16 \\
  $H^7$ & 2 & 1 & 8 & 10 & 21 & 19 & 25 & 30 & 23 & 34 & 43 & 30 \\
  \hline
  \, & $\phi_{35}^{31}$ & $\phi_{35}^{4}$ & $\phi_{35}^{22}$ & $\phi_{35}^{13}$ & $\phi_{56}^{30}$ & $\phi_{56}^{3}$ & $\phi_{70}^{18}$ & $\phi_{70}^{9}$ & $\phi_{84}^{12}$ & $\phi_{84}^{15}$ & $\phi_{105}^{5}$ & $\phi_{105}^{26}$ \\
  $H^0$ &$\cdot$&$\cdot$&$\cdot$&$\cdot$&$\cdot$&$\cdot$&$\cdot$&$\cdot$&$\cdot$&$\cdot$&$\cdot$&$\cdot$\\
  $H^1$ &$\cdot$& 1 &$\cdot$&$\cdot$&$\cdot$&$\cdot$&$\cdot$&$\cdot$&$\cdot$&$\cdot$&$\cdot$&$\cdot$\\
  $H^2$ &$\cdot$& 2 &$\cdot$&$\cdot$&$\cdot$& 1 &$\cdot$&$\cdot$&$\cdot$&$\cdot$& 1 &$\cdot$\\
  $H^3$ &$\cdot$& 3 &$\cdot$& 1 &$\cdot$& 3 &$\cdot$& 1 & 1 &$\cdot$& 4 &$\cdot$\\
  $H^4$ &$\cdot$& 9 & 2 & 4 &$\cdot$& 9 & 5 & 6 & 9 & 5 & 15 & 1 \\
  $H^5$ & 3 & 30 & 16 & 14 & 11 & 30 & 30 & 30 & 50 & 27 & 53 & 29 \\
  $H^6$ & 23 & 63 & 45 & 36 & 53 & 70 & 86 & 80 & 127 & 78 & 122 & 113 \\
  $H^7$ & 43 & 52 & 47 & 44 & 74 & 71 & 101 & 85 & 117 & 108 & 134 & 137 \\
  \hline
  \, & $\phi_{105}^{12}$ & $\phi_{105}^{6}$ & $\phi_{105}^{15}$ & $\phi_{105}^{21}$ & $\phi_{120}^{4}$ & $\phi_{120}^{25}$ & $\phi_{168}^{6}$ & $\phi_{168}^{21}$ & $\phi_{189}^{22}$ & $\phi_{189}^{20}$ & $\phi_{189}^{5}$ & $\phi_{189}^{7}$ \\
  \hline
  $H^0$ &$\cdot$&$\cdot$&$\cdot$&$\cdot$&$\cdot$&$\cdot$&$\cdot$&$\cdot$&$\cdot$&$\cdot$&$\cdot$&$\cdot$\\
  $H^1$ &$\cdot$&$\cdot$&$\cdot$&$\cdot$&$\cdot$&$\cdot$&$\cdot$&$\cdot$&$\cdot$&$\cdot$&$\cdot$&$\cdot$\\
  $H^2$ &$\cdot$& 1 &$\cdot$&$\cdot$& 2 &$\cdot$& 1 &$\cdot$&$\cdot$&$\cdot$& 1 &$\cdot$\\
  $H^3$ & 2 & 7 &$\cdot$&$\cdot$& 9 &$\cdot$& 7 &$\cdot$&$\cdot$&$\cdot$& 6 & 4 \\
  $H^4$ & 14 & 27 & 5 & 1 & 33 &$\cdot$& 36 & 2 & 5 & 7 & 25 & 23 \\
  $H^5$ & 63 & 78 & 34 & 20 & 99 & 19 & 128 & 35 & 61 & 73 & 90 & 86 \\
  $H^6$ & 154 & 160 & 101 & 94 & 194 & 99 & 267 & 145 & 215 & 233 & 216 & 205 \\
  $H^7$ & 147 & 159 & 133 & 124 & 185 & 136 & 249 & 200 & 255 & 251 & 239 & 243 \\
  \hline
  \, & $\phi_{189}^{17}$ & $\phi_{189}^{10}$ & $\phi_{210}^{10}$ & $\phi_{210}^{6}$ & $\phi_{210}^{13}$ & $\phi_{210}^{21}$ & $\phi_{216}^{16}$ & $\phi_{216}^{9}$ & $\phi_{280}^{8}$ & $\phi_{280}^{17}$ & $\phi_{280}^{18}$ & $\phi_{280}^{9}$ \\
  \hline
  $H^0$ &$\cdot$&$\cdot$&$\cdot$&$\cdot$&$\cdot$&$\cdot$&$\cdot$&$\cdot$&$\cdot$&$\cdot$&$\cdot$&$\cdot$\\
  $H^1$ &$\cdot$&$\cdot$&$\cdot$&$\cdot$&$\cdot$&$\cdot$&$\cdot$&$\cdot$&$\cdot$&$\cdot$&$\cdot$&$\cdot$\\
  $H^2$ &$\cdot$&$\cdot$&$\cdot$& 2 &$\cdot$&$\cdot$&$\cdot$&$\cdot$& 2 &$\cdot$&$\cdot$&$\cdot$\\
  $H^3$ &$\cdot$& 5 & 4 & 13 &$\cdot$&$\cdot$& 1 & 3 & 9 &$\cdot$&$\cdot$& 4 \\
  $H^4$ & 6 & 33 & 32 & 51 & 9 & 2 & 13 & 21 & 47 & 7 & 19 & 28 \\
  $H^5$ & 55 & 125 & 128 & 157 & 68 & 45 & 99 & 87 & 191 & 73 & 126 & 121 \\
  $H^6$ & 182 & 277 & 295 & 326 & 214 & 185 & 287 & 224 & 427 & 257 & 351 & 306 \\
  $H^7$ & 226 & 276 & 307 & 313 & 253 & 248 & 296 & 275 & 404 & 343 & 388 & 345 \\
  \hline
\, & $\phi_{315}^{16}$ & $\phi_{315}^{7}$ & $\phi_{336}^{14}$ & $\phi_{336}^{11}$ & $\phi_{378}^{14}$ & $\phi_{378}^{9}$ & $\phi_{405}^{15}$ & $\phi_{405}^{8}$ & $\phi_{420}^{10}$ & $\phi_{420}^{13}$ & $\phi_{512}^{12}$ & $\phi_{512}^{11}$ \\ 
  \hline
  $H^0$ &$\cdot$&$\cdot$&$\cdot$&$\cdot$&$\cdot$&$\cdot$&$\cdot$&$\cdot$&$\cdot$&$\cdot$&$\cdot$&$\cdot$\\
  $H^1$ &$\cdot$&$\cdot$&$\cdot$&$\cdot$&$\cdot$&$\cdot$&$\cdot$&$\cdot$&$\cdot$&$\cdot$&$\cdot$&$\cdot$\\
  $H^2$ &$\cdot$&$\cdot$&$\cdot$&$\cdot$&$\cdot$&$\cdot$&$\cdot$&$\cdot$&$\cdot$&$\cdot$&$\cdot$&$\cdot$\\
  $H^3$ &$\cdot$& 4 & 2 & 2 & 2 & 2 &$\cdot$& 11 & 7 &$\cdot$& 6 & 2 \\
  $H^4$ & 21 & 31 & 33 & 26 & 33 & 28 & 12 & 73 & 61 & 20 & 61 & 29 \\
  $H^5$ & 141 & 136 & 179 & 128 & 188 & 150 & 118 & 268 & 258 & 139 & 290 & 180 \\
  $H^6$ & 393 & 347 & 456 & 344 & 498 & 405 & 397 & 588 & 598 & 417 & 710 & 524 \\
  $H^7$ & 441 & 385 & 468 & 416 & 533 & 467 & 480 & 598 & 602 & 510 & 731 & 627
 \end{tabular}
 \end{table}
\end{center}  
 
 \subsection{The root system $E_8$}
 
  Let $\Phi$ be the root system of type $E_8$. Then the Poincar\'e polynomial
 of $X_{\Phi}$ is
 \begin{equation*}
 \begin{array}{rl}
   P(X_{\Phi},t)= &  1313187309t^8 + 818120000t^7 + 235845616t^6 + 37527168t^5 + \\
   &+3539578t^4 + 202496t^3 + 6888t^2 + 128t + 1.
 \end{array}
 \end{equation*}
 The cohomology of $X_{\Phi}$ as a representation of the Weyl group
 of $E_8$ is given in Table~\ref{e8cohtable}.
 
\begin{center}
\begin{longtable}[ht!]{rrrrrrrrr}
\caption{The cohomology groups of the complement of the toric arrangement associated to $E_8$ as
 representations of the Weyl group.}\\ 
 \label{e8cohtable} 

\, & $\phi_{1}^{0}$ & $\phi_{1}^{120}$ & $\phi_{8}^{91}$ & $\phi_{8}^{1}$ 
   & $\phi_{28}^{8}$ & $\phi_{28}^{68}$ & $\phi_{35}^{74}$ & $\phi_{35}^{2}$ \\*
\nobreakhline
$H^0$ & 1 &$\cdot$&$\cdot$&$\cdot$&$\cdot$&$\cdot$&$\cdot$&$\cdot$\\*
$H^1$ & 1 &$\cdot$& 1 &$\cdot$&$\cdot$&$\cdot$&$\cdot$& 1 \\*
$H^2$ &$\cdot$&$\cdot$& 1 &$\cdot$&$\cdot$& 1 &$\cdot$& 2 \\*
$H^3$ &$\cdot$&$\cdot$&$\cdot$&$\cdot$&$\cdot$& 2 &$\cdot$& 2 \\*
$H^4$ &$\cdot$&$\cdot$&$\cdot$&$\cdot$&$\cdot$& 3 &$\cdot$& 3 \\*
$H^5$ &$\cdot$&$\cdot$&$\cdot$&$\cdot$&$\cdot$& 6 &$\cdot$& 9 \\*
$H^6$ &$\cdot$&$\cdot$& 1 &$\cdot$& 2 & 19 & 1 & 31 \\*
$H^7$ & 2 &$\cdot$& 8 & 2 & 22 & 55 & 22 & 72 \\*
$H^8$ & 3 & 1 & 15 & 10 & 48 & 68 & 56 & 83 \\
\hline
\pagebreak[3] \, & $\phi_{50}^{56}$ & $\phi_{50}^{8}$ & $\phi_{56}^{49}$ & $\phi_{56}^{19}$
   & $\phi_{70}^{32}$ & $\phi_{84}^{64}$ & $\phi_{84}^{4}$ & $\phi_{112}^{63}$ \\*
\nobreakhline
$H^0$ &$\cdot$&$\cdot$&$\cdot$&$\cdot$&$\cdot$&$\cdot$&$\cdot$&$\cdot$\\*
$H^1$ &$\cdot$&$\cdot$&$\cdot$&$\cdot$&$\cdot$&$\cdot$& 1 &$\cdot$\\*
$H^2$ &$\cdot$&$\cdot$&$\cdot$&$\cdot$&$\cdot$&$\cdot$& 2 & 1 \\*
$H^3$ &$\cdot$& 1 & 1 &$\cdot$&$\cdot$&$\cdot$& 3 & 2 \\*
$H^4$ &$\cdot$& 4 & 2 &$\cdot$& 1 &$\cdot$& 7 & 3 \\*
$H^5$ &$\cdot$& 13 & 4 & 1 & 6 &$\cdot$& 21 & 9 \\*
$H^6$ & 4 & 37 & 14 & 7 & 29 & 5 & 70 & 40 \\* 
$H^7$ & 47 & 92 & 52 & 39 & 101 & 66 & 166 & 123 \\* 
$H^8$ & 93 & 115 & 97 & 89 & 147 & 145 & 196 & 202 \\
\hline
\pagebreak[3] \, & $\phi_{112}^{3}$ & $\phi_{160}^{7}$ & $\phi_{160}^{55}$ & $\phi_{168}^{24}$ 
   & $\phi_{175}^{36}$ & $\phi_{175}^{12}$ & $\phi_{210}^{52}$ & $\phi_{210}^{4}$ \\*
\nobreakhline
$H^0$ &$\cdot$&$\cdot$&$\cdot$&$\cdot$&$\cdot$&$\cdot$&$\cdot$&$\cdot$\\*
$H^1$ &$\cdot$&$\cdot$&$\cdot$&$\cdot$&$\cdot$&$\cdot$&$\cdot$&$\cdot$\\*
$H^2$ &$\cdot$&$\cdot$& 1 &$\cdot$&$\cdot$&$\cdot$&$\cdot$& 1 \\*
$H^3$ &$\cdot$&$\cdot$& 2 &$\cdot$&$\cdot$& 1 &$\cdot$& 4 \\*
$H^4$ &$\cdot$&$\cdot$& 4 & 1 &$\cdot$& 7 &$\cdot$& 13 \\*
$H^5$ &$\cdot$&$\cdot$& 13 & 11 & 3 & 27 &$\cdot$& 47 \\*
$H^6$ & 3 & 10 & 52 & 72 & 44 & 105 & 28 & 160 \\*
$H^7$ & 58 & 96 & 166 & 248 & 214 & 293 & 198 & 395 \\*
$H^8$ & 167 & 246 & 284 & 354 & 348 & 384 & 380 & 482 \\
\hline
 \, & $\phi_{300}^{44}$ & $\phi_{300}^{8}$ & $\phi_{350}^{38}$ & $\phi_{350}^{14}$ & 
     $\phi_{400}^{7}$ & $\phi_{400}^{43}$ & $\phi_{420}^{20}$ & $\phi_{448}^{25}$ \\*
\nobreakhline
$H^0$ &$\cdot$&$\cdot$&$\cdot$&$\cdot$&$\cdot$&$\cdot$&$\cdot$&$\cdot$\\*
$H^1$ &$\cdot$&$\cdot$&$\cdot$&$\cdot$&$\cdot$&$\cdot$&$\cdot$&$\cdot$\\*
$H^2$ &$\cdot$&$\cdot$&$\cdot$&$\cdot$&$\cdot$&$\cdot$&$\cdot$&$\cdot$\\*
$H^3$ & 2 &$\cdot$&$\cdot$& 2 &$\cdot$& 1 &$\cdot$&$\cdot$\\*
$H^4$ & 12 &$\cdot$&$\cdot$& 11 &$\cdot$& 8 & 3 & 2 \\*
$H^5$ & 54 & 3 & 9 & 51 &$\cdot$& 34 & 32 & 15 \\*
$H^6$ & 206 & 59 & 102 & 214 & 33 & 130 & 184 & 98 \\*
$H^7$ & 536 & 324 & 425 & 587 & 276 & 402 & 618 & 392 \\*
$H^8$ & 674 & 568 & 682 & 765 & 643 & 699 & 883 & 755 \\

\hline
\, & $\phi_{448}^{39}$ & $\phi_{448}^{9}$ & $\phi_{525}^{12}$ & $\phi_{525}^{36}$ & 
     $\phi_{560}^{5}$ & $\phi_{560}^{47}$ & $\phi_{567}^{6}$ & $\phi_{567}^{46}$ \\*
\nobreakhline
$H^0$ &$\cdot$&$\cdot$&$\cdot$&$\cdot$&$\cdot$&$\cdot$&$\cdot$&$\cdot$\\*
$H^1$ &$\cdot$&$\cdot$&$\cdot$&$\cdot$&$\cdot$&$\cdot$&$\cdot$&$\cdot$\\*
$H^2$ &$\cdot$&$\cdot$&$\cdot$&$\cdot$&$\cdot$& 1 & 2 &$\cdot$\\*
$H^3$ & 1 &$\cdot$& 3 &$\cdot$&$\cdot$& 4 & 9 &$\cdot$\\*
$H^4$ & 5 &$\cdot$& 17 &$\cdot$&$\cdot$& 11 & 29 &$\cdot$\\*
$H^5$ & 29 & 1 & 77 & 11 &$\cdot$& 44 & 112 & 3 \\*
$H^6$ & 144 & 48 & 323 & 152 & 41 & 189 & 407 & 106 \\*
$H^7$ & 452 & 326 & 883 & 642 & 363 & 581 & 1028 & 588 \\*
$H^8$ & 781 & 727 & 1148 & 1026 & 882 & 988 & 1278 & 1052\\
\hline
\, & $\phi_{700}^{42}$ & $\phi_{700}^{6}$ & $\phi_{700}^{16}$ & $\phi_{700}^{28}$ & 
     $\phi_{840}^{13}$ & $\phi_{840}^{31}$ & $\phi_{840}^{26}$ & $\phi_{840}^{14}$ \\*
\nobreakhline
$H^0$ &$\cdot$&$\cdot$&$\cdot$&$\cdot$&$\cdot$&$\cdot$&$\cdot$&$\cdot$\\*
$H^1$ &$\cdot$&$\cdot$&$\cdot$&$\cdot$&$\cdot$&$\cdot$&$\cdot$&$\cdot$\\*
$H^2$ &$\cdot$& 1 &$\cdot$&$\cdot$&$\cdot$&$\cdot$&$\cdot$&$\cdot$\\*
$H^3$ &$\cdot$& 6 & 1 &$\cdot$&$\cdot$&$\cdot$&$\cdot$&$\cdot$\\*
$H^4$ &$\cdot$& 29 & 11 & 1 &$\cdot$& 5 & 1 & 8 \\*
$H^5$ & 5 & 132 & 76 & 26 & 9 & 44 & 37 & 88 \\*
$H^6$ & 140 & 494 & 373 & 251 & 132 & 239 & 311 & 461 \\*
$H^7$ & 759 & 1246 & 1103 & 948 & 669 & 809 & 1134 & 1331 \\*
$H^8$ & 1324 & 1560 & 1496 & 1422 & 1386 & 1449 & 1699 & 1790 \\
\hline
\, & $\phi_{972}^{32}$ & $\phi_{972}^{12}$ & $\phi_{1008}^{39}$ & $\phi_{1008}^{9}$ & 
     $\phi_{1050}^{34}$ & $\phi_{1050}^{10}$ & $\phi_{1134}^{20}$ & $\phi_{1296}^{33}$ \\*
\nobreakhline
$H^0$ &$\cdot$&$\cdot$&$\cdot$&$\cdot$&$\cdot$&$\cdot$&$\cdot$&$\cdot$\\*
$H^1$ &$\cdot$&$\cdot$&$\cdot$&$\cdot$&$\cdot$&$\cdot$&$\cdot$&$\cdot$\\*
$H^2$ &$\cdot$&$\cdot$&$\cdot$&$\cdot$&$\cdot$& 1 &$\cdot$&$\cdot$\\*
$H^3$ &$\cdot$& 2 & 3 &$\cdot$&$\cdot$& 6 &$\cdot$& 2 \\*
$H^4$ & 1 & 20 & 15 &$\cdot$&$\cdot$& 31 & 7 & 14 \\*
$H^5$ & 26 & 127 & 70 & 4 & 20 & 158 & 82 & 76 \\*
$H^6$ & 293 & 574 & 316 & 122 & 293 & 654 & 507 & 368 \\*
$H^7$ & 1229 & 1604 & 1007 & 733 & 1284 & 1760 & 1667 & 1244 \\*
$H^8$ & 1933 & 2111 & 1757 & 1623 & 2061 & 2288 & 2369 & 2236 \\
\hline
\, & $\phi_{1296}^{13}$ & $\phi_{1344}^{8}$ & $\phi_{1344}^{38}$ & $\phi_{1344}^{19}$ & 
     $\phi_{1400}^{11}$ & $\phi_{1400}^{29}$ & $\phi_{1400}^{32}$ & $\phi_{1400}^{8}$ \\*
\nobreakhline
$H^0$ &$\cdot$&$\cdot$&$\cdot$&$\cdot$&$\cdot$&$\cdot$&$\cdot$&$\cdot$\\*
$H^1$ &$\cdot$&$\cdot$&$\cdot$&$\cdot$&$\cdot$&$\cdot$&$\cdot$&$\cdot$\\*
$H^2$ &$\cdot$& 2 &$\cdot$&$\cdot$&$\cdot$&$\cdot$&$\cdot$&$\cdot$\\*
$H^3$ &$\cdot$& 9 &$\cdot$&$\cdot$& 1 &$\cdot$&$\cdot$& 6 \\*
$H^4$ &$\cdot$& 44 &$\cdot$& 1 & 11 &$\cdot$&$\cdot$& 44 \\*
$H^5$ & 15 & 222 & 19 & 36 & 77 & 12 & 32 & 213 \\*
$H^6$ & 203 & 889 & 325 & 300 & 393 & 211 & 391 & 857 \\*
$H^7$ & 1016 & 2329 & 1549 & 1203 & 1339 & 1122 & 1710 & 2349 \\*
$H^8$ & 2124 & 2969 & 2587 & 2282 & 2413 & 2323 & 2751 & 3067 \\
\hline
\, & $\phi_{1400}^{37}$ & $\phi_{1400}^{7}$ & $\phi_{1400}^{20}$ & $\phi_{1575}^{10}$ & 
     $\phi_{1575}^{34}$ & $\phi_{1680}^{22}$ & $\phi_{2016}^{19}$ & $\phi_{2100}^{28}$ \\*
\nobreakhline
$H^0$ &$\cdot$&$\cdot$&$\cdot$&$\cdot$&$\cdot$&$\cdot$&$\cdot$&$\cdot$\\*
$H^1$ &$\cdot$&$\cdot$&$\cdot$&$\cdot$&$\cdot$&$\cdot$&$\cdot$&$\cdot$\\*
$H^2$ &$\cdot$&$\cdot$&$\cdot$&$\cdot$&$\cdot$&$\cdot$&$\cdot$&$\cdot$\\*
$H^3$ & 3 &$\cdot$&$\cdot$&$\cdot$& 7 &$\cdot$&$\cdot$&$\cdot$\\*
$H^4$ & 19 &$\cdot$& 10 &$\cdot$& 47 & 12 & 3 & 4 \\*
$H^5$ & 98 & 5 & 110 & 38 & 236 & 132 & 56 & 101 \\*
$H^6$ & 441 & 163 & 634 & 450 & 971 & 768 & 448 & 780 \\*
$H^7$ & 1398 & 1019 & 2048 & 1920 & 2645 & 2454 & 1798 & 2822 \\*
$H^8$ & 2439 & 2261 & 2914 & 3083 & 3445 & 3486 & 3419 & 4239 \\
\hline
\, & $\phi_{2100}^{16}$ & $\phi_{2100}^{20}$ & $\phi_{2240}^{28}$ & $\phi_{2240}^{10}$ & 
     $\phi_{2268}^{10}$ & $\phi_{2268}^{30}$ & $\phi_{2400}^{17}$ & $\phi_{2400}^{23}$ \\*
\nobreakhline
$H^0$ &$\cdot$&$\cdot$&$\cdot$&$\cdot$&$\cdot$&$\cdot$&$\cdot$&$\cdot$\\*
$H^1$ &$\cdot$&$\cdot$&$\cdot$&$\cdot$&$\cdot$&$\cdot$&$\cdot$&$\cdot$\\*
$H^2$ &$\cdot$&$\cdot$&$\cdot$&$\cdot$&$\cdot$&$\cdot$&$\cdot$&$\cdot$\\*
$H^3$ & 3 &$\cdot$&$\cdot$& 4 & 5 &$\cdot$&$\cdot$&$\cdot$\\*
$H^4$ & 35 & 13 & 1 & 43 & 50 &$\cdot$& 3 & 10 \\*
$H^5$ & 236 & 160 & 72 & 288 & 311 & 71 & 59 & 95 \\*
$H^6$ & 1129 & 958 & 720 & 1303 & 1357 & 720 & 494 & 587 \\*
$H^7$ & 3304 & 3076 & 2882 & 3654 & 3733 & 2863 & 2060 & 2183 \\*
$H^8$ & 4479 & 4365 & 4473 & 4840 & 4910 & 4482 & 4022 & 4081 \\
\hline
\, & $\phi_{2688}^{20}$ & $\phi_{2800}^{25}$ & $\phi_{2800}^{13}$ & $\phi_{2835}^{14}$ & 
     $\phi_{2835}^{22}$ & $\phi_{3150}^{18}$ & $\phi_{3200}^{22}$ & $\phi_{3200}^{16}$ \\*
\nobreakhline
$H^0$ &$\cdot$&$\cdot$&$\cdot$&$\cdot$&$\cdot$&$\cdot$&$\cdot$&$\cdot$\\*
$H^1$ &$\cdot$&$\cdot$&$\cdot$&$\cdot$&$\cdot$&$\cdot$&$\cdot$&$\cdot$\\*
$H^2$ &$\cdot$&$\cdot$&$\cdot$&$\cdot$&$\cdot$&$\cdot$&$\cdot$&$\cdot$\\*
$H^3$ & 1 &$\cdot$&$\cdot$& 2 &$\cdot$&$\cdot$&$\cdot$& 1 \\*
$H^4$ & 18 & 15 &$\cdot$& 32 & 7 & 20 & 12 & 26 \\*
$H^5$ & 197 & 134 & 48 & 278 & 150 & 232 & 189 & 283 \\*
$H^6$ & 1205 & 743 & 519 & 1453 & 1121 & 1420 & 1314 & 1593 \\*
$H^7$ & 3941 & 2616 & 2326 & 4369 & 3941 & 4628 & 4505 & 4881 \\*
$H^8$ & 5604 & 4792 & 4655 & 5999 & 5798 & 6570 & 6568 & 6746 \\
\hline
\, & $\phi_{3240}^{9}$ & $\phi_{3240}^{31}$ & $\phi_{3360}^{25}$ & $\phi_{3360}^{13}$ & $\phi_{4096}^{12}$ & $\phi_{4096}^{26}$ & $\phi_{4096}^{27}$ & $\phi_{4096}^{11}$ \\*
\nobreakhline
$H^0$ &$\cdot$&$\cdot$&$\cdot$&$\cdot$&$\cdot$&$\cdot$&$\cdot$&$\cdot$\\*
$H^1$ &$\cdot$&$\cdot$&$\cdot$&$\cdot$&$\cdot$&$\cdot$&$\cdot$&$\cdot$\\*
$H^2$ &$\cdot$&$\cdot$&$\cdot$&$\cdot$&$\cdot$&$\cdot$&$\cdot$&$\cdot$\\*
$H^3$ &$\cdot$& 3 & 1 &$\cdot$& 6 &$\cdot$& 2 &$\cdot$\\*
$H^4$ &$\cdot$& 30 & 19 & 1 & 70 & 5 & 25 &$\cdot$\\*
$H^5$ & 23 & 193 & 157 & 53 & 488 & 163 & 203 & 55 \\*
$H^6$ & 461 & 953 & 885 & 606 & 2293 & 1424 & 1124 & 698 \\*
$H^7$ & 2517 & 3152 & 3141 & 2794 & 6568 & 5390 & 3888 & 3337 \\*
$H^8$ & 5319 & 5605 & 5755 & 5600 & 8795 & 8220 & 7040 & 6790 \\
\hline
\, & $\phi_{4200}^{24}$ & $\phi_{4200}^{12}$ & $\phi_{4200}^{15}$ & $\phi_{4200}^{21}$ & $\phi_{4200}^{18}$ & $\phi_{4480}^{16}$ & $\phi_{4536}^{23}$ & $\phi_{4536}^{13}$ \\*
\nobreakhline
$H^0$ &$\cdot$&$\cdot$&$\cdot$&$\cdot$&$\cdot$&$\cdot$&$\cdot$&$\cdot$\\*
$H^1$ &$\cdot$&$\cdot$&$\cdot$&$\cdot$&$\cdot$&$\cdot$&$\cdot$&$\cdot$\\*
$H^2$ &$\cdot$&$\cdot$&$\cdot$&$\cdot$&$\cdot$&$\cdot$&$\cdot$&$\cdot$\\*
$H^3$ &$\cdot$& 5 &$\cdot$&$\cdot$& 1 &$\cdot$&$\cdot$&$\cdot$\\*
$H^4$ & 8 & 64 & 2 & 12 & 24 & 28 & 16 & 1 \\*
$H^5$ & 187 & 463 & 87 & 155 & 306 & 337 & 188 & 77 \\*
$H^6$ & 1532 & 2255 & 850 & 1028 & 1908 & 2026 & 1175 & 845 \\*
$H^7$ & 5663 & 6622 & 3630 & 3843 & 6168 & 6574 & 4234 & 3825 \\*
$H^8$ & 8510 & 8971 & 7065 & 7158 & 8743 & 9337 & 7767 & 7592 \\
\hline
\, & $\phi_{4536}^{18}$ & $\phi_{5600}^{19}$ & $\phi_{5600}^{15}$ & $\phi_{5600}^{21}$ & $\phi_{5670}^{18}$ & $\phi_{6075}^{14}$ & $\phi_{6075}^{22}$ & $\phi_{7168}^{17}$ \\*
\nobreakhline
$H^0$ &$\cdot$&$\cdot$&$\cdot$&$\cdot$&$\cdot$&$\cdot$&$\cdot$&$\cdot$\\
$H^1$ &$\cdot$&$\cdot$&$\cdot$&$\cdot$&$\cdot$&$\cdot$&$\cdot$&$\cdot$\\*
$H^2$ &$\cdot$&$\cdot$&$\cdot$&$\cdot$&$\cdot$&$\cdot$&$\cdot$&$\cdot$\\*
$H^3$ &$\cdot$&$\cdot$&$\cdot$&$\cdot$&$\cdot$&$\cdot$& 4 &$\cdot$\\*
$H^4$ & 30 & 9 & 4 & 17 & 36 & 15 & 69 & 11 \\*
$H^5$ & 349 & 167 & 129 & 212 & 431 & 329 & 597 & 208 \\*
$H^6$ & 2070 & 1261 & 1145 & 1382 & 2580 & 2421 & 3113 & 1605 \\*
$H^7$ & 6643 & 4973 & 4817 & 5120 & 8309 & 8435 & 9358 & 6376 \\*
$H^8$ & 9428 & 9470 & 9397 & 9533 & 11794 & 12403 & 12852 & 12136 
\end{longtable}
\end{center}

\subsection*{Acknowledgements}
This paper is based on parts of my thesis \cite{bergvallthesis} written 
at Stockholms Universitet.
The author would like to thank Carel Faber and Jonas Bergström for helpful discussions
and comments. The author would also like to thank Alessandro Oneto, Ivan Martino and Travis Scrimshaw for useful comments on early versions of this manuscript
and Federico Ardila for pointing out a preprint of the interesting paper \cite{ardilacastillohenley}.
Thanks also goes to Emanuele Delucchi for interesting discussions.
Some of the computations for the $E_8$-case were performed on resources provided by SNIC
through PDC under projects PDC-2016-4 and PDC-2018-74.
\clearpage
\bibliographystyle{amsalpha}

\renewcommand{\bibname}{References} 

\clearpage
\bibliography{references} 
\end{document}